\documentclass[12pt,a4paper]{article}
\usepackage{amsmath,amsfonts,amssymb,amsthm,amscd,latexsym}

\newcommand{\map}{\rightarrow}

\newcommand{\ul}[1]{\underline{#1}}

\newcommand{\cs}{\text{c}} 

\def\CO{\mathcal{O}}
\def\CQ{\mathcal{Q}}
\def\CE{\mathcal{E}}
\def\CA{\mathcal{A}}
\def\CC{\mathcal{C}}
\def\CH{\mathcal{H}}

\def\dbar{\bar{\partial}}

\def\rank{\mathop{\mathrm{Rank}}\nolimits}
\def\ker{\mathop{\mathrm{Ker}}\nolimits}
\def\coker{\mathop{\mathrm{Coker}}\nolimits}
\def\im{\mathop{\mathrm{Im}}\nolimits}

\newtheorem{theorem}{Theorem}[section]
\newtheorem{lemma}[theorem]{Lemma}
\newtheorem{proposition}[theorem]{Proposition}
\newtheorem{definition}[theorem]{Definition}

\author{Toby N. Bailey and Liana David}

\title{The Penrose transform for compactly supported cohomology.}

\begin{document}
\maketitle

\renewcommand{\thefootnote}{}
\footnotetext{Mathematics Subject Classification. Primary 32L25;
  Secondary 22E46, 53C28.}  

\begin{abstract}
Let the manifold $X$ parametrise a family of compact complex
submanifolds of the complex (or CR) manifold $Z$.   Under mild
conditions the Penrose transform typically provides isomorphisms 
between a cohomology group of a holomorphic vector bundle $V \map Z$
and the kernel of a differential operator between sections of vector
bundles over $X$. When the spaces in question are homogeneous for a
group $G$ the Penrose transform provides an intertwining operator
between representations.   

In this paper we develop a Penrose transform  for compactly supported
cohomology on $Z$. We provide a number of examples where a compactly
supported cohomology group is shown to be isomorphic to the cokernel
of a differential operator between compactly supported sections of 
vector bundles over $X$.   We consider also how the ``Serre duality'' 
pairing carries through the transform.  
\end{abstract}

\section{Introduction} 
Let the real manifold $X$ parametrise a family of compact
holomorphic  submanifolds of a complex (or perhaps CR) manifold $Z$. 
If certain conditions are satisfied, the {\em Penrose transform} 
enables one to interpret Dolbeault cohomology with values in a
holomorphic vector bundle on $Z$ in terms of kernels and cokernels of
differential operators on $X$.  This procedure was developed
in a representation-theoretic context by Schmid \cite{S}
and independently by Penrose  \cite{P} as part of his ``Twistor Programme''. 
(A cohomological interpretation of Penrose's work 
did not appear until \cite{EPW}, which considers the
transform only in the holomorphic category.)

Let the set of points in $X$ incident with each $z \in Z$ be
contractible and of dimension $d$. Let $V \map Z$ be a holomorphic
vector bundle.  Then under conditions which appear below there is a
spectral sequence
 \[
    E^{p,q}_1 = \Gamma(X, V_{p,q}) \Longrightarrow 
       H^{p+q}(Z, V).  
 \]     
where the $V_{p,q}$ are certain complex vector bundles over $X$ and 
the maps at the $E_1$ level are first order differential operators. 
Our main result is that in these circumstances one also has a spectral
sequence 
  \[
    E^{p,q}_1 = \Gamma_c(X, V_{p,q}) \Longrightarrow 
       H^{p+q-d}_c(Z, V).  
   \]   
where the ``c'' subscripts refer to compactly supported sections and
cohomology. 

We also consider the ``Serre duality'' 
pairing between cohomology and  compactly supported cohomology
and the resulting pairing between the kernels and cokernels of
differential operators arising from the transform.

\section{Involutive structures and cohomology} \label{is}

We recall briefly the results about involutive structures and their
cohomology that we will need.  Details can be found in \cite{BES}. 
We then define the compactly supported involutive cohomology. 
Our conventions are that $TM, \CE_M^k $ will always refer to the {\em
  complexified} tangent bundle and the $k$-th exterior power of the
complexified  cotangent bundles of the real manifold $M$.  In
particular, $\CE_M = \CE^0$ refers to the trivial complex line bundle 
whose sections are complex-valued functions on $M$.  We omit 
the ``$M$'' and write simply $\CE^k$, etc, if there can be no
confusion as to which manifold is intended.

\begin{definition}
  An \ul{involutive structure} on the smooth manifold $M$ is a complex
  sub-bundle $T^{0,1}\subset TM$ such that $[T^{0,1},T^{0,1}]\subset
  T^{0,1}$ (meaning that the space of smooth sections is closed under
  Lie bracket).  Define the vector bundle $\mathcal E^{1,0}
  \subseteq \CE^1$ to be the annihilator of $T^{0,1}$ and define
  $\CE^{0,1}$ by the exactness of 
    $$
     0\rightarrow\mathcal E^{1,0}\rightarrow\mathcal
    E^{1} \rightarrow\mathcal E^{0,1}\rightarrow 0.
    $$
We write $\mathcal E^{p,q} =  \wedge^{p}\mathcal E^{1,0} \otimes
 \wedge^{q}   \mathcal E^{0,1} $.
\end{definition}

\begin{definition}
A complex vector bundle $V\rightarrow M$ is \ul{compatible} with the
involutive structure $\mathcal E$ (or $\mathcal E$-compatible) 
if there is defined a linear operator 
$$
\dbar :\Gamma (M,V)\rightarrow \Gamma (M,V\otimes\mathcal E^{0,1})
$$
such that 
$$
\dbar (fs)=f\dbar (s)+(\dbar f)s,\quad\forall f\in\mathcal E
(M),s\in\Gamma (M,V)
$$
and such that the extension to 
 $$
    \dbar :
 \Gamma (M,V\otimes\mathcal E^{0,q})
  \rightarrow \Gamma (M,V\otimes\mathcal E^{0,q+1})
$$
satisfies $\dbar^2=0$. 
\end{definition}

\begin{definition}
Given an $\mathcal E$-compatible vector bundle $V\rightarrow M$ 
define the \ul{involutive cohomology} $H^{*}_{\mathcal E}(M,V)$ to be the 
cohomology of the complex
$$
   \Gamma (M,V)\stackrel{\dbar}{\rightarrow}\Gamma (M,\mathcal
E^{0,1}\otimes V)\stackrel{\dbar}{\rightarrow}\Gamma (M,\mathcal
E^{0,2}\otimes V)\stackrel{\dbar}{\rightarrow}\cdots
$$ 
Writing $\Gamma_{\cs}(M,V)$ for the compactly supported smooth sections of
$V$, etc, we define the  \ul{compactly supported 
involutive cohomology} $H^{*}_{\mathcal E,\cs}(M,V)$ to be the 
cohomology of the complex
$$
   \Gamma_{\cs} (M,V)\stackrel{\dbar}{\rightarrow}\Gamma_{\cs} (M,\mathcal
E^{0,1}\otimes V)\stackrel{\dbar}{\rightarrow}\Gamma_{\cs} (M,\mathcal
E^{0,2}\otimes V)\stackrel{\dbar}{\rightarrow}\cdots
$$ 
\end{definition}
We will be using the following examples of involutive structures: 
\begin{enumerate} 
\item   An involutive structure $T^{0,1}\subset TM$ is a {\em 
 complex structure}  iff
$TM=T^{0,1}\oplus\overline{T^{0,1}}$. In this case the compatible vector
bundles are the holomorphic ones and the involutive
cohomology is the Dolbeault cohomology. The compactly supported 
cohomology is also the usual compactly supported Dolbeault cohomology. 
\item 
If $T^{0,1}\cap \overline{T^{0,1}} = \{ 0 \}$ then the involutive
structure is a {\em CR-structure}. This is the structure acquired by a
real hypersurface in a complex manifold.  Our definition is wider than
most since it includes ``higher codimension'' cases and includes also
complex manifolds. 
The involutive cohomology is often known in this case as 
$\dbar_{b}$ cohomology (the ``b'' standing for boundary).  
\item 
Let $\eta :F\rightarrow Z$ be a fibre bundle. Then $\mathcal
C^{1,0}=\eta^{*}(\mathcal E^{1}_Z)$ is an involutive structure on
$F$. Bundles on $F$ which are pull-backs of  bundles on $M$ are
compatible.    Suppose the fibres of $\eta$ have finite-dimensional de
Rham cohomology. Then the $k$-th cohomology of the fibre defines a
vector bundle $\CH^k \map Z$.     The involutive cohomology in this
case (often called ``relative de Rham cohomology'') is  
fiber-wise de Rham
cohomology, parametrised by $Z$.    To be precise, we have 
 \[
     H^p_\CC(F, \eta^* V) = \Gamma(Z, V \otimes \CH^p).
 \]
\end{enumerate} 

The principal calculational tool of \cite{BES} 
concerns the situation where one has
two involutive structures $\mathcal A^{1,0}\subset \mathcal E^{1,0}$
on a manifold $M$.  Defining 
$B^{1}:=\mathcal E^{1,0}/\mathcal A^{1,0}$ we also have a short exact
sequence 
 \[
    0 \map B^1 \map \CA^{0,1} \map \CE^{0,1} \map 0.  
 \]
Using this to filter the complex for computing $H^*_\CA(M,V)$ we
obtain:  
\begin{proposition} \label{ss} 
Let $\mathcal A^{1,0}\subset \mathcal E^{1,0}$ be two involutive
structures on the manifold $M$ and let $V\rightarrow M$ be an $\mathcal
A$-compatible vector bundle. Then $V$ is also $\mathcal E$-compatible,
 and there is a spectral sequence
$$
E_{1}^{p,q}=H^{q}_{\mathcal E}(M,B^{p}\otimes V)\Longrightarrow
H^{p+q}_{\mathcal A}(M,V)
$$
where 
$B^{p}:=\wedge^{p}B^{1}$.

Taking instead the compactly supported cohomology the same filtration
gives a spectral sequence 
 $$
E_{1}^{p,q}=H^{q}_{\mathcal E, c}(M,B^{p}\otimes V)\Longrightarrow
H^{p+q}_{\mathcal A, c}(M,V)
$$
\end{proposition}
The proof is straightforward  (see  \cite{BES} for details) 
and the compactly
supported case is completely analogous. 

\section{The Penrose transform} \label{pt} 

We summarise the Penrose transform here as presented in \cite{BES} 
in a simple case. (See also \cite{BD,CC}.) 
The initial data are a double fibration of smooth oriented manifolds 

  \begin{equation}  \label{cor} 
  \begin{picture}(100,70)
\put(50,60){\makebox(0,0){$F$}}
\put(40,50){\vector(-3,-4){20}}
\put(60,50){\vector(3,-4){20}}
\put(15,15){\makebox(0,0){$Z$}}
\put(85,15){\makebox(0,0){$X$}}
\put(25,40){\makebox(0,0){$\scriptstyle \eta$}}
\put(75,40){\makebox(0,0){$\scriptstyle \tau$}}
\end{picture}
\end{equation}
with the following properties.
\begin{enumerate} 
\item 
$Z$ has an involutive structure $\CQ$ which is a CR-manifold. 
(We recall that our definition includes the
possibility that $Z$ is complex.) 
  \item 
The maps $\eta$ and $\tau$ are fiber-bundle projections such that 
$\tau$ has compact complex fibers.  
\item 
The map $\eta$ embeds the fibres of $\tau$ as holomorphic submanifolds of
$Z$. (A submanifold of a CR-manifold is holomorphic if the involutive
structure of $Z$ restricts to give a complex structure on the
submanifold.)  
\end{enumerate} 
In this situation, we can endow $F$ with the involutive structure
$\CE$ defined by $\CE^{1,0}$ being the annihilator of the $\tau$-vertical  
vectors that are of type $(0,1)$ with respect to the complex structure
on the fibres of $\tau$.

Let $V\rightarrow Z$ be a $\mathcal Q$-compatible vector bundle
on $Z$.     The Penrose transform proceeds in three stages.

\paragraph{Step 1}
Define an involutive structure $\CA$ on $F$ by $\CA^{1,0} = \eta^*
\CQ^{1,0}$.   Then $\CQ$-cohomology on $Z$ and $\CA$-cohomology on $F$
can be related as follows. 
Introduce first the 
involutive structure $\CC^{1,0} = \eta^* \CE^1_Z$ on $F$. 
This is exactly the third example in \S\ref{is}.   On $F$ we have 
 \[
       0 \map \CA^{1,0} \map \CC^{1,0} \map \eta^*\CQ^{0,1} \map 0 
 \]
and the standard spectral sequence of Proposition~\ref{ss} gives 
 \[
   E_1^{p,q} = H^q_\CC (F, \eta^* V \otimes \eta^* \CQ^{0,p}) 
   \Longrightarrow H^{p+q}_\CA(F, \eta^* V).  
 \]
Using the identification of $H_\CC$ in \S\ref{is} we obtain 
 \[
     E_2^{p,q} = H^p_\CQ(Z, V\otimes \CH^q). 
 \]
A particular case of importance is where the fibres of $\eta$ are
contractible and we deduce that 
 \[
          H^p_\CQ(Z, V)  \cong  H^p_\CA(Z, V) 
 \]
with the map being given by pull-back of a representative form. 

\paragraph{Step 2}
We define $B^1 \map F$ by $B^{1}:=\mathcal E^{1,0}/\mathcal A^{1,0}$
and employ the standard spectral sequence of Proposition~\ref{ss} 
to obtain 
 \[
    E^{p,q}_1 = H^q_\CE(F, \eta^*V\otimes B^p) \Longrightarrow 
       H^{p+q}_\CA (F,\eta^* V).  
 \]        

\paragraph{Step 3}
We identify the $\CE$-cohomology groups appearing in the previous
step.   The $\CE$-cohomology is the Dolbeault cohomology of the
fibres of $\tau$ parametrised over $X$.   Let us assume that 
a complex $\CE$-compatible vector bundle $E\map F$ is such that the
dimension of the Dolbeault cohomology $H^p(\tau^{-1}(x), \left. E
\right|_{\tau^{-1}x})$ is constant as $x\in X$ varies.  Then this
cohomology defines a vector bundle $\tau_*^p E \map X$. 
(It is in fact the $p$-th 
direct image of the sheaf of smooth sections of $E$ holomorphic on
each fibre of $\tau$.)    Then 
 \[
     H^p_{\CE}(F, E) = \Gamma(X, \tau^p_*E).
 \]

Defining 
 \[
       V_{p,q} = \tau^q_* (\eta^* V \otimes B^p) 
 \]
we thus have 
    \[
    E^{p,q}_1 = \Gamma(X, V_{p,q}) \Longrightarrow 
       H^{p+q}_\CA (F,\eta^* V).  
 \]     
The maps are first order differential operators.

Combining the above steps {\em in the case where $\eta$ has
  contractible fibres} {we arrive at the {\em Penrose transform} 
which is the spectral sequence 
 \[
    E^{p,q}_1 = \Gamma(X, V_{p,q}) \Longrightarrow 
       H^{p+q}_\CQ (Z, V).  
 \]

\section{The compactly supported transform}

We retain  the setting of \S\ref{pt},  supposing also that the fibers of
$\eta$ have  finite dimensional compactly supported de Rham cohomology.
We proceed by analogy with the real Penrose transform, using the fact
that the spectral sequence of Proposition~\ref{ss} is valid also for
compactly supported cohomology. 
In the first step we will need to identify the compactly supported
cohomology of the involutive structure $\CC$ on $F$. 
\begin{lemma}
Consider the involutive 
structure $\mathcal C^{1,0}=\eta^{*}(\mathcal E^{1}_{Z})$ on $F$ 
(as in Step 1 for the real Penrose transform). 
 For $V\map Z$ a vector bundle and $ k\geq 0$,
\[
  H^{k}_{\mathcal C,\cs}(F, \eta^*V)
 \cong\Gamma_{c}(Z,\mathcal H^{k}_{c}\otimes V)
\] 
where $\mathcal H^{k}_{c}$ is the bundle whose fiber over $z\in Z$ is the
$k$-compactly supported de Rham cohomology of $\eta^{-1}(z)$. 
\end{lemma}

\begin{proof}
We will use only the case where the fibres of $\eta$ are contractible
in this paper, when one can   apply fiber by fiber the
homotopy formula for compactly supported cohomology (see e.g.\ 
\cite[\S 4]{BT}). (The general case then follows from a 
``Cech de Rham complex'' argument.) 
\end{proof}

\paragraph{Step 1}
Following Step 1 of the Penrose transform in the compactly supported
case and using the above Lemma to identify the $\CC$-cohomology we
obtain the pull-back spectral sequence: 
 \[
 E_{2}^{p,q}=H^{p}_{\mathcal Q,\cs}(Z,\mathcal H^{q}_{c}\otimes V)   
  \Longrightarrow
  H^{p+q}_{\mathcal A,\cs}(F,\eta^{*}V).
 \]

If the fibers of $\eta$ are contractible this immediately converges 
and we get 
 \[
  H^{k}_{\mathcal Q,\cs}(Z,V)\cong H^{k+d}_{\mathcal A ,\cs}(F,\eta^{*}V)
 \]
where $d$ is the dimension of the fibres of $\eta$. 
Choose  
$\rho \in \Gamma(F, \CC^{0,d})$ which represents ``$1$''
in the compactly supported top-degree cohomology of each fibre. 
The map is given in this case by pull-back of forms followed by
wedging with $\rho$.

\paragraph{Step 2}
We use the standard spectral sequence of Proposition~\ref{ss}
for compactly
supported cohomology  to obtain  
 $$
 E_{1}^{p,q}=H^{q}_{\mathcal E 
 ,\cs}(F,B^{p}\otimes\eta^{*}V)\Longrightarrow H^{p+q}_{\mathcal
 A ,\cs}(F,\eta^{*}V)
 $$

\paragraph{Step 3}
We identify 
$H^{q}_{\mathcal E ,\cs}(F,\eta^{*}V\otimes B^{p})
  \cong\Gamma_{c}(X,V_{p,q})$ where $V_{p,q}$ are exactly the same
  bundles as arise in the standard transform.  
This follows almost immediately from
  the corresponding fact for the usual transform. 
When the fibres of $\eta$ are contractible, combining the steps proves
the following. 
\begin{theorem} 
In the situation of \S\ref{pt} in the case 
where the fibres of $\eta$ are 
contractible there is a spectral sequence 
 $$
  E_{1}^{p,q}=\Gamma_{c}(X,V_{p,q})
    \Longrightarrow H^{p+q-d}_{\mathcal Q ,\cs}(Z,V).
 $$
The vector bundles that appear are exactly those for the standard
transform and the differential operators in the spectral sequence are
the same as those that appear in the non-compactly supported case,
but acting between  compactly supported sections of the relevant
bundles.
\end{theorem} 

\section{The  bilinear pairing}

For an involutive structure $\CA$, define 
$\kappa_\CA$ to be the 
top exterior power of $\CA^{1,0}$ (by analogy with the definition of
the holomorphic canonical bundle on a complex manifold).

\begin{definition} 
On a manifold $F$  with involutive structure 
$\CA$ and compatible vector bundle $V$, let $k+l = \rank \CA^{0,1}$. 
  The \underline{natural bilinear pairing} between involutive cohomology and
  compactly supported involutive cohomology on $F$ 
 \[
\int _{F} : H^{k}_{\mathcal A ,c}(F,V) \times H^{l}_{\mathcal A}(F,
  V^{*} \otimes \kappa_\CA )\rightarrow\mathbb{C}, 
  \]
is that given by wedge product of representative forms (combined with
contraction between the vector space $V$ and its dual $V^*$) 
followed by integration. 
\end{definition} 
When  $F$ is a complex manifold and  $\mathcal Q$ is the 
complex structure this is the Serre duality pairing 
that, when $F$ is compact,  
identifies the (necessarily finite-dimensional) cohomology spaces
as mutually dual.   

Let $F$ now be the correspondence space for the Penrose transform as
previously. We have also on $F$ the corresponding pairing for
$\CE$-cohomology and a compatible vector bundle $E$: 
   $$
   \int _{F}:H^{q}_{\mathcal E ,c}(F, E)\times H^{r}_{\mathcal
   E}(F,  E^*\otimes \kappa_{\mathcal E})\rightarrow
   \mathbb{C}.
  $$
whenever $q+r=\rank (\mathcal \CE^{0,1})$ (which is the complex
dimension of the fibres of $\tau$). 

\begin{lemma}
Let  $\kappa_\tau$ be the line bundle on $F$ which restricts to each
fibre of $\tau$ to be the canonical bundle of that fibre. (Recall that
the fibres of $\tau$ are naturally complex manifolds.)  Then 
 \[
      \kappa_\CE = \tau^* ( \Lambda^{\rm{top}}_X ) \otimes \kappa_\tau. 
 \]
(Here $\Lambda^{\mathrm{top}}_X$ denotes the line bundle of
complex-valued top-degree forms on $X$.) 
\end{lemma}
\begin{proof}
In the circumstances of the Penrose transform as we have been
discussing,  there is a short exact sequence of vector bundles 
on $F$ 
 \[
 0 \map \tau^* \CE^1_X \map \CE^{1,0} \map \CE^{1,0}_\tau \map 0 
 \]
where $\CE^{1,0}_\tau$ denotes the vector bundle of forms of type 
$(1,0)$ in the complex structure of the fibres of $\tau$ (so that 
$\CE^{1,0}_\tau = \overline{\CE^{0,1}}$).  The result
follows by taking top exterior powers. 
\end{proof}

\begin{proposition} 
Consider the situation in the Penrose transform where we have 
a $\CE$-compatible vector bundle $E \map F$ such that 
 \[
 H^q_\CE(F, E) = \Gamma(X,\tau_*^q E), \quad 
  H^r_{\CE,c}(F, E^*\otimes \kappa_\CE ) 
   = \Gamma_c(X,\tau_*^r (E^*\otimes \kappa_\CE ) ). 
 \]
When $q+r= \rank (\CE^{0,1})$ we can identify 
 \[
 \tau_*^r (E^*\otimes \kappa_\CE ) = (\tau_*^q E)^* \otimes      
 \Lambda^{\rm{top}}_X 
 \]
and the pairing on $F$ for $\CE$ cohomology is given by 
 \[
     \int_X : \Gamma_c(X,\tau_*^q E) \times 
          \Gamma(X,  (\tau_*^q E)^* \otimes \Lambda^{\rm{top}}_X) 
   \map \Bbb C, 
 \]
where we are contracting the vector bundle $\tau_*^q E$ with its dual
and integrating the resulting top-order form over $X$. 
\end{proposition} 

\begin{proof} 
Note that by the preceding Lemma, 
 \[
   E^* \otimes \kappa_\CE = E^* \otimes \kappa_\tau \otimes
   \tau^* \Lambda^{\text{top}}_X 
 \]
and so 
 \[
 \tau_*^r (E^* \otimes \kappa_\CE) = \tau_*^r (E^* \otimes
\kappa_\tau) \otimes  \Lambda^{\text{top}}_X.
 \] 
For each fibre $\tau^{-1}x$ the cohomology groups 
 \[
     H^q(\tau^{-1}x, E) \quad \text{and} \quad  H^r(\tau^{-1}x,
     E^*\otimes \kappa_\tau)
 \]
are Serre-dual and so we can identify 
 \[
        \tau_*^r (E^* \otimes \kappa_\CE ) = (\tau_*^q E)^*.  
 \]
Now split the pairing integral on $F$ into a fibre integral, which is
precisely the Serre duality pairing, followed by an integral over the
base. 
\end{proof}

\begin{proposition}
Let 
 $$
  E_{1}^{p,q}=\Gamma_{c}(X,V_{p,q})\Longrightarrow 
  H^{p+q}_{\mathcal A, c}(F,\eta^{*}(V))
 $$
    and  
 $$
   \widetilde{E_{1}}^{s,t}=\Gamma (X,(V^{*}\otimes
  \kappa_{\mathcal Q})_{s,t})   \Longrightarrow 
   H^{s+t}_{\mathcal  A}(F,\eta^{*}(V^{*}\otimes \kappa_{\CA}))
 $$
be the spectral sequences for the Penrose transform as discussed
above. There is a pairing for $r\geq 1$  
 \[
    E_r^{p,q} \times   \widetilde{E_{1}}^{s,t} \map \Bbb C,  
 \quad    p+s = \rank(B^1), \, q+t = \rank(\CE^{0,1}) 
 \]
which for $r=1$ is the pairing for $\CE$-cohomology on $F$ and which
converges to the $\CA$-cohomology pairing on $F$. 
\end{proposition}

\begin{proof}
Note first that $\eta^* \kappa_\CQ = \kappa_\CA$. 
The proof uses the fact that the filtrations of the spectral sequences
are induced by a sub-bundle of the $\mathcal A^{0,1}$. One can check
directly that $\int_{F}$ descends to the corresponding terms of
the $r$-level of the spectral sequences.
\end{proof}

A similar analysis  holds for the pull-back stage of our double
fibration where 
instead of the involutive structure $\mathcal E^{1,0}$ we consider the
involutive structure $\mathcal C^{1,0}$.      The outcome in the case
where the fibres of $\eta$ are contractible is that the pairing 
 \[
    H^{k-d}_{\CQ,c}(Z,V) \times H^l_\CQ(Z, V^*\otimes \kappa_\CQ) \map \Bbb C
 \]
(where $k+l-d = \rank \CQ^{0,1}$ and $d$ is the dimension of the 
fibres of $\tau$) pulls back to give the pairing 
 \[
    H^{k}_{\CA,c}(F,\eta^*(V)) \times 
  H^l_\CA(F, \eta^*(V^*)\otimes \kappa_\CA) \map \Bbb C. 
 \]

\section{Examples}

\subsection{Euclidean space $\mathbb{R}^{3}$}
Let $X= \mathbb{R}^{3}$ and let $Z$ be the  
total space of the holomorphic tangent bundle of $\Bbb CP_1$, thought
of as the parametrisation space of oriented straight lines in $X$. 
One obtains a double fibration where $F$ is the space of ``points on
oriented lines in $\Bbb R^3$''. 

Let $G$ be the double cover of the group of Euclidean motions of
$X=\Bbb R^3$. We realise $G$ as  
 \[
  G = \{  (A,B) \mid  A \in SU(2), \text{$B$ is a $2\times 2$ trace-free 
   hermitian matrix}     \}
 \]
acting on $X= \Bbb R^3 =$ the space of trace-free Hermitian $2\times 2$
matrices $x$ according to $x \mapsto AxA^* + B$.

As a homogeneous space $Z= G/L$ where 
  \[
    L=  \left\{ \left( 
      \begin{pmatrix}
        e^{i\theta}  & 0 \\ 
        0  & e^{-i\theta}
      \end{pmatrix} , 
 \begin{pmatrix}
        z  & 0 \\ 
        0  &  -z 
      \end{pmatrix} \right) ,  \theta\in \Bbb R, z \in \Bbb C 
     \right\}.
\]
For $n \in \Bbb Z, \lambda \in \Bbb C$ define 
$\CO(n,\lambda)$ to be  the $G$-homogeneous holomorphic line
bundle on $Z$ associated to the character $e^{-in\theta-\lambda z}$ of
$L$.  Taking $\CQ$ to be the complex structure on $Z$, so that the
involutive cohomology is the usual Dolbeault cohomology,  it is easy to
check that $\kappa_\CQ = \CO(-4,0)$.

The Penrose transform gives  isomorphisms \cite{CC}  
\begin{gather*}
     H^1_\CQ(Z,\CO(-2,\lambda)) \overset{=}{\longrightarrow} 
    \ker( \Delta + 2\lambda^2)  \\ 
    H^2_\CQ(Z,\CO(-2,\lambda)) \overset{=}{\longrightarrow} 
    \coker( \Delta + 2\lambda^2)   
\end{gather*}
where $\Delta$ denotes the Laplacian mapping the space of  smooth functions on
$\mathbb{R}^{3}$ to itself. 

The fibres of $\eta$ in this case are 1-dimensional, and so we
immediately obtain the following. 
\begin{theorem}
Let $\Delta_{c}$ denote the mapping from the space of compactly
supported smooth functions to itself given by the Laplacian. 
Then 
\begin{gather*}
     H^0_{\CQ,c}(Z,\CO(-2,\lambda)) \overset{=}{\longrightarrow} 
    \ker( \Delta_c + 2\lambda^2)  \\ 
    H^1_\CQ(Z,\CO(-2,\lambda)) \overset{=}{\longrightarrow} 
    \coker( \Delta_c + 2\lambda^2)  
\end{gather*}
(The first observation is trivial since both sides are zero.) 
\end{theorem} 

The ``Serre duality pairing'' on $Z$ 
 \[
   H^1_{\CQ,c}(Z,\CO(-2,\lambda)) \times 
    H^1_{\CQ}(Z,\CO(-2, - \lambda))   \map \Bbb C 
 \]
translates into the pairing 
 \[
     (f, [g] ) \mapsto \int_{\Bbb R^3} fg 
 \]
where $f \in \ker( \Delta + 2\lambda^2)$ and 
$g \in \coker( \Delta_c + 2\lambda^2)$.

\subsection{A CR example}

This is Penrose's original transform, restricted to real Minkowski
space (see \cite{nmjw}). 

On $C^4$ with coordinates $z_1, \dots , z_4$ define 
the pseudo-hermitian form $\Phi$ by 
 \[
 \Phi(z,z) = z_1 \bar{z}_3 + z_2 \bar{z}_4
   + z_3 \bar{z}_1 + z_4 \bar{z}_2.
 \]
Let $I$ denote the projective line $z_3=z_4=0$ and let $Z$ be the
5-dimensional CR-manifold 
 \[
    Z = \{ z \in \Bbb CP_3 \, | \,  \Phi(z,z)=0 \} \setminus I .
 \]
The space of complex projective lines which lie in $Z$ can be
identified with ``Minkowski space'' $X = \Bbb R^4$ with Lorentzian
metric. 

Taking $\CQ$ to be the CR-structure on $Z$, there are $\CQ$-compatible line
bundles $\CO(n), n \in \Bbb Z$ which are the restrictions of the usual
holomorphic line bundles on $\Bbb CP_3$. We have $\kappa_\CQ =
\CO(-4)$.  The fibres of $\eta$ are again 1-dimensional and so we have
a ``dimension shift'' of one for the compactly supported cohomology. 

\subsubsection{The case of $\CO(-2)$} 

The Penrose transform gives an isomorphism 
 \[
     H^1_\CQ(Z, \CO(-2)) \overset{=}{\longrightarrow} 
            \ker \Box 
 \]
where $\Box$ is the wave operator associated to the Lorentz structure
on $\Bbb R^4$ mapping the space of smooth functions to itself.

We deduce immediately that the compactly supported transform gives 
isomorphisms 
\begin{gather*}
  H^0_{\CQ,c}(Z, \CO(-2))  \overset{=}{\longrightarrow}  \ker \Box_c
  \\ 
 H^1_{\CQ,c}(Z, \CO(-2))  \overset{=}{\longrightarrow}  \coker \Box_c
\end{gather*}
where $\Box_c$ denotes the wave operator  
mapping the space of compactly supported  smooth functions to itself. 
The first isomorphism is trivial, both sides being zero. 

The ``Serre duality pairing'' on $Z$ 
 \[
   H^1_{\CQ,c}(Z,\CO(-2)) \times 
    H^1_{\CQ}(Z,\CO(-2)) \map \Bbb C  
 \]
translates into the pairing 
 \[
     (f, [g] ) \mapsto \int_{\Bbb R^4} fg 
 \]
where $f \in \ker \Box$ and 
$g \in \coker \Box_c$. 

\subsubsection{The case of $\CO(-1)$ and $\CO(-3)$} 

We will not consider this case in detail.  The Penrose transform 
gives 
 \[
   H^1_{\CQ}(Z, \CO(-3))  \overset{=}{\longrightarrow}  \ker D^- 
 \]
where $D^-$ is the Dirac operator from smooth sections of the spin
bundle $S^-$ to the other spin bundle $S^+$.   

The compactly supported transform gives 
   \[
   H^1_{\CQ,c}(Z, \CO(-1))  \overset{=}{\longrightarrow}  \coker D^+_c
 \]
where $D^+_c$ is the Dirac operator from compactly supported 
smooth sections of the spin
bundle $S^+ $ to compactly supported smooth sections of $S^-$.

The Serre duality pairing between these groups becomes the pairing 
 \[
      (\alpha, [\beta]) \mapsto \int_{\Bbb R^4} \epsilon(\alpha,
      \beta) 
 \]
where $\alpha \in \ker D^-$ and $[\beta] \in \coker D^+$ and $\epsilon$
is the complex bilinear skew form on $S^-$. 
   
\subsection{The case of $\CO(-4)$ and $\CO(0)$} 

We recall that the smooth complex 2-forms $\Lambda^2$ on $X=\Bbb R^4$
split as a direct sum of self-dual and anti-self-dual: 
 \[
    \Lambda^2 = \Lambda^2_+ \oplus \Lambda^2_- .
 \]
 The Penrose transform gives 
 \[
   H^1_{\CQ}(Z, \CO(-4))  \overset{=}{\longrightarrow}  \ker 
 (d: \Lambda^2_- \map \Lambda^3 ) .  
 \]

The Penrose transform on $\CO(0)$ gives 
\begin{gather*}
    H^1_{\CQ}(Z, \CO(0))  \overset{=}{\longrightarrow} 
  \frac{ \ker  (d: \Lambda^1 \map \Lambda^2_- )} 
  {\im d : \Lambda^0 \map \Lambda^1}  \\ 
    H^2_{\CQ}(Z, \CO(0))  \overset{=}{\longrightarrow} 
  \frac{ \Lambda^2_-} 
  {\im d : \Lambda^1 \map \Lambda^2_-} 
\end{gather*}
and so we can immediately deduce that 
 \begin{gather*}
    H^0_{\CQ,c}(Z, \CO(0))  \overset{=}{\longrightarrow} 
  \frac{ \ker  (d: \Lambda^1_{c} \map \Lambda^2_{-,c} )} 
  {\im d : \Lambda^0_c \map \Lambda^1_c}  \\ 
    H^1_{\CQ,c}(Z, \CO(0))  \overset{=}{\longrightarrow} 
  \frac{ \Lambda^2_{-,c}} 
  {\im d : \Lambda^1_c \map \Lambda^2_{-,c}} 
\end{gather*}
The cohomology group on the left of the first statement is clearly
zero and so we deduce that a compactly supported complex-valued 1-form
with the self-dual part of its exterior derivative vanishing is
necessarily the exterior derivative of a compactly supported
function. 

The ``Serre duality pairing'' on $Z$ 
 \[
   H^1_{\CQ,c}(Z,\CO(0)) \times 
    H^1_{\CQ}(Z,\CO(-4)) \map \Bbb C  
 \]
translates into the pairing 
 \[
     (f, [g] ) \mapsto \int_{\Bbb R^4} \langle f, g \rangle  
 \]
where $f, g$ are anti-self-dual 2-forms and $\langle f, g \rangle $ 
is the usual bilinear form.

\subsection{Odd-dimensional hyperbolic spaces}

In \cite{BD} a twistor correspondence and Penrose transform
for $X =$ hyperbolic space of dimension $2n+1$ is described. It is
equivariant with respect to $G= \mathrm{Spin}_0(2n+1,1)$.  The space
$Z$ is an open orbit in the isotropic Grassmanian of complex
$n$-planes in $\Bbb C^{2n+2}$.  There are $G$-homogeneous holomorphic line
bundles $\CO(n,\lambda), n \in \Bbb Z, \lambda \in \Bbb C$ on $Z$.

The smooth Penrose transform gives isomorphisms between Dolbeault
cohomologies and kernels and cokernels of operators on $H^{2n+1}$. 
(We omit the ``$\CQ$'' for the involutive structure which is the
complex structure on $Z$.) 

\begin{gather*}
  H^{\frac{n(n+1)}{2}}(Z ,\mathcal O(-2n,\lambda ))   
      \overset{=}{\longrightarrow}  
       \ker(\Delta -(\lambda^{2}-n^{2}))\\
  H^{\frac{n(n+1)}{2}+1}(Z,\mathcal O(-2n,\lambda ))  
       \overset{=}{\longrightarrow} 
       \coker(\Delta -(\lambda^{2}-n^{2}))\\
  H^{n}(Z,\mathcal O(-2,\lambda ))   
        \overset{=}{\longrightarrow}  
         \ker(\Delta -(\lambda^{2}-n^{2}))\\
   H^{n+1}(Z,\mathcal O(-2,\lambda))   
     \overset{=}{\longrightarrow}  
          \coker(\Delta -(\lambda^{2}-n^{2}))
\end{gather*}
where $\Delta$ denotes the hyperbolic Laplacian defined from the space
of smooth functions 
on $H^{2n+1}$ to itself.

Similar isomorphisms hold for the compactly supported Dolbeault
cohomologies, with a shift in dimension by $1$ (which is again the
dimension of the fibres of $\eta$): 
  \begin{gather*}
      H^{\frac{n(n+1)}{2}-1}_{c}(Z,\CO(-2n,\lambda)) 
              \overset{=}{\longrightarrow}  
              \ker(\Delta_{c} -(\lambda^{2}-n^{2}))\\
     H^{\frac{n(n+1)}{2}}_{c}(Z,\CO(-2n,\lambda ))   
              \overset{=}{\longrightarrow} 
                    \coker(\Delta_{c}-(\lambda^{2}-n^{2}))\\
    H^{n-1}_{c}(Z,\CO(-2,\lambda))   
             \overset{=}{\longrightarrow}  
              \ker(\Delta_{c}-(\lambda^{2}-n^{2}))\\
    H^{n}_{c}(Z,\CO(-2,\lambda ))  
                 \overset{=}{\longrightarrow}  
            \coker(\Delta_{c}-(\lambda^{2}-n^{2}))
\end{gather*}
where $\Delta_{c}$ denotes the hyperbolic Laplacian from the
compactly supported smooth functions to itself. 

In this case, $\kappa_q = \CO(-2n-2,0)$ on $Z$ which has complex
dimension $n(n+3)/2$ and so we have Serre-duality pairings 
\begin{gather*}
  H^{n}_{c}(Z,\CO(-2,\lambda ))   \times
   H^{\frac{n(n+1)}{2}}(Z,\CO(-2n,-\lambda ))
   \rightarrow\mathbb{C}  \\
 H^{\frac{n(n+1)}{2}}_c(Z,\CO (-2n, \lambda))
\times
 H^{n}(Z,\CO(-2, - \lambda ))\times
\rightarrow\mathbb{C}
\end{gather*}
which both become the pairing 
 \[
  \coker (\Delta_c -(\lambda^{2}-n^{2}))\times
\ker (\Delta -(\lambda^{2}-n^{2}))  \rightarrow\mathbb{C} 
 \]
induced by multiplication of functions and then integration.  The
situation with the cohomology  giving rise to eigenspinors of the
Dirac operator in \cite{BD} is very similar.

\end{document}